\documentclass{conm-p-l}
\usepackage{bbm}
\usepackage{amssymb,amsthm,amsfonts,amsmath,float}
\newtheorem{theorem}{Theorem}[section]

\theoremstyle{definition}
\newtheorem{definition}[theorem]{Definition}

\theoremstyle{remark}
\newtheorem{remark}[theorem]{Remark}

\numberwithin{equation}{section}
\newcommand\reals{{\mathbb R}}

\newcommand\nc\newcommand
\nc\cC{{\mathcal C}}
\nc{\remove}[1]{}

\begin{document}

\title{New bounds for equiangular lines}

\author{Alexander Barg}
\address{Dept. of ECE and Inst. for Systems Research, University of Maryland, College Park, MD, 20742, and Inst. for Problems of Inform. Trans., RAS, Moscow, Russia}
\email{abarg@umd.edu}

\author{Wei-Hsuan Yu}
\address{Department of Math. and Inst. for Systems Research, University of Maryland, College Park, MD 20742}
\email{mathyu@math.umd.edu}

\date{\today}
\dedicatory{To Ilya Dumer, on the occasion of his 60th birthday}

\begin{abstract}
A set of lines in $\reals^n$ is called equiangular if the angle
between each pair of lines is the same. We address the question of
determining the maximum size of equiangular line sets in
$\mathbb{R}^n,$ using semidefinite programming to improve the
upper bounds on this quantity. Improvements are obtained in
dimensions $24 \leq n \leq 136$. In particular, we show that the
maximum number of equiangular lines in $\reals^n$ is $276$ for all
$24 \leq n \leq 41$  and is 344 for $n=43.$ This provides a
partial resolution of the conjecture set forth by Lemmens and
Seidel (1973).
\end{abstract}

\maketitle

\section{Introduction}
A set of lines in a metric space is called equiangular if the
angle between each pair of lines is the same. We are interested in
upper bounds on the number of equiangular lines in $\mathbb{R}^n$.
In other words, if we have a set of unit vectors $S=
\{x_i\}_{i=1}^M$ and there is a constant $c>0$ such that $|\langle
x_i, x_j \rangle| = c$ for all $ 1 \leq i \neq j \leq M$, what is
the maximum cardinality of $S$? Denote this quantity by $M(n).$ The problem of determining
$M(n)$ looks elementary but
a general answer has so far proved elusive: Until recently the maximum number of
equiangular lines in $\mathbb{R}^n$ was known only for $16$ values
of the dimension $n$. The history of this problem started with
Hanntjes \cite{Hann48} who found $M(n)$ for $n=2$ and $3$ in $1948$.
Van Lint and Seidel \cite{sei66} found the largest number
of equiangular lines for $4\le n \leq 7$. In $1973$, Lemmens and
Seidel \cite{lem73} used linear-algebraic methods to determine $M(n)$ for most values of $n$ in the region
$8\le n\le 23$.
Gerzon (see \cite{lem73}) gave the following upper on $M(n)$.
\begin{theorem}[Gerzon]
If there are $M$ equiangular lines in $\mathbb{R}^n$, then
  \begin{equation}\label{eq:gerzon}
M \leq \frac{n(n+1)}{2}
  \end{equation}
\end{theorem}
Gerzon's upper bound can be attained only for a very small number
of values of $n.$ Currently, such constructions are known only for
$n=2,3,7,$ and $23.$ Neumann (see \cite{lem73}, Theorem 3.2)
proved a fundamental result in this area:

\begin{theorem}[Neumann]\label{thm:Newmann}
If there are $M$ equiangular lines in $\mathbb{R}^n$ with angle
$\arccos \alpha$ and $M>2n$, then $1/\alpha$ is an odd integer.
\end{theorem}
Note that if $M$ attains the Gerzon bound, then $(n+2)\alpha^2=1$
\cite[Thm.3.5]{lem73}. Therefore, if the cardinality of an
equiangular line set attains the Gerzon bound, then $n$ has to be
$2$ or $3$ or an odd square minus two and the angle between pairs
of lines is $\arccos {1}/({\sqrt{n+2}})$.

A set of unit vectors $S=\{x_1,x_2,\dots\}\subset \reals^n$ is called {\em two-distance} if
$\langle x_i,x_j\rangle\in\{a,b\}$ for some $a,b$ and all $i\ne j.$

\begin{theorem}[Larman, Rogers, and Seidel \cite{lar77}]\label{thm:LRS}
Let $S$ be a spherical two-distance set in $\mathbb{R}^n.$ If $|S| >
2n+3$ and $a> b$, then $b=\frac{ka-1}{k-1}$ for some integer $k$ such that $2 \leq k \leq
(1+\sqrt{2n})/{2}$.
\end{theorem}
The condition $|S|>2n+3$ was improved to $|S|>2n+1$ by Neumaier
\cite{neu81}. He also gave an example of a two-distance set with
cardinality $2n+1$ that violates the integeraity condition of $k$.
This example is obtained from the spherical embedding of the
conference graph.

If the spherical two-distance set gives rise to equiangular lines, then $a=-b,$ so Theorem
\ref{thm:LRS} implies that $a=1/(2k-1),$ which is the statement of the Neumann theorem.
The assumption of Theorem \ref{thm:LRS} is more restrictive than of Theorem \ref{thm:Newmann}, but
in return we obtain an upper bound on $k$. For instance, if
$n=40$, then $k$ can be only $2$ or $3$, so the angle has to be
$\arccos \alpha,$ where $\alpha={1/3}$ or $1/5.$ The assumption of Theorem
\ref{thm:LRS} is satisfied since there exist equiangular line sets with $M\ge 2n+4$
for all $n\ge 15.$



\begin{table}[t]\begin{center}{ \begin{tabular}{|c|c|c|c|c|c|c|}
\hline
 $n$ & $M(n)$ & $1/\alpha$ && $n$ &$M(n)$ & $1/\alpha$ \\
\hline
2&3&2& &17&  48-50 &5 \\
3&6&$\sqrt{5}$& &18 &  48-61 &5\\
4&6&3; $\sqrt 5$&& 19&72-76&5 \\
5&10&3&& 20 &90-96 &5\\
6&16&3& &21 &126&5\\
$7\le n\le 13$&28&3&&22&176&5\\
14&28-30&3;\,5  &&$23$&276&5 \\
15&36&5 && $24 \leq n \leq 42$& $\geq$ 276 &5\\
16& 40-42 &5 & & 43 & $\geq$ 344&7 \\
\hline
\end{tabular}}
\end{center}\caption{Known bounds on $M(n)$ in small dimensions}\label{table:known}
\end{table}
\begin{table}[t]
\begin{center}\begin{tabular}{|c|c|c|c|c|c|c|}
\hline
 $n$ & $M(n)$ & SDP bound  & &$n$ &$M(n)$ & SDP bound \\
\hline
3&6& 6& &18 & 48-61 &61\\
4&6& 6&& 19&72-76&76 \\
5&10&10&& 20 &90-96 &96\\
6&16&16& &21 &126&126\\
$7\le n\le 13$&28&28&&22&176&176\\
14&28-30&30 &&$23  $&276&276 \\
15&36&36    && $24 \leq n \leq 41$ & 276&276\\
16& 40-42&42 && 42&$\geq$ 276&288 \\
17& 48-50&51 && 43& 344& 344 \\
 \hline
\end{tabular}
\end{center}\caption{Bounds on $M(n)$ including new results}\label{table:new}
\vspace*{-.2in}\end{table}

The known bounds on $M(n)$ for small dimensions are summarized
in Table \ref{table:known} \cite{lem73}, \cite{sus07} (the latter for the upper bound on $M(17)$); in particular, $M(n)$ was known exactly only if
$2\le n\le13; n= 15,21,22,23$. In the
unsettled cases the best known upper bound on $M(n)$ is usually
the Gerzon bound. Lemmens and Seidel \cite[Thm.~4.5]{lem73} further showed that
  \begin{equation}\label{eq:1/3}
  M_{1/3}(n) \leq 2(n-1), \quad n \geq 16,
  \end{equation}
 where $M_{\alpha}(n)$ is the maximum size of an equiangular line set
when the value of the angle is $\arccos \alpha.$
They also conjectured that $M_{1/5}(n)= 276$ for $23 \leq n \leq 185,$
observing that if this conjecture is true, then $M(n)=276$ for $24 \leq n \leq 41$ and $M(43)=344.$
Note that generally we have \cite{lem73}:
  \begin{equation}\label{eq:relative}
    M_\alpha(n)\le\frac {n(1-\alpha^2)}{1-n\alpha^2}
  \end{equation}
 valid for all $\alpha$ such that the denominator is positive. This inequality is sometimes called the {\em relative bound} as opposed to
  the ``absolute bound'' of \eqref{eq:gerzon}.

In this paper we use the semidefinite programming (SDP) method to derive some new bounds
on $M(n).$ Our main results are summarized in Table \ref{table:new}.
In particular, exact values of $M(n)$ are obtained for $24\le n\le 41$
and for $n=43$ where previous results gave divergent bounds:
we show that $M(n)=276$ for $24 \leq n \leq 41$ and $M(43)=344$.
These results are established by performing
computations with SDP. We also show that $M_{1/5}(n) = 276$ for $23 \leq n
\leq 60$. These results resolve a part of the Lemmens-Seidel conjecture and enable us to obtain the results in
Table \ref{table:new}. For $44 \leq n \leq 136$,
we also obtain new upper bounds on $M(n),$ improving upon the Gerzon bound, although no new
exact values are found in this range.
Below in the paper we give a more complete table of the computation results.

An interesting question relates to the asymptotic behavior of $M(n)$ for $n\to \infty.$
For a long time the best known constructions were able to attain the growth order of $M(n)=\Omega(n),$
until D. de Caen \cite{Caen00} constructed a family of $\frac{2}{9}(n+1)^2$
equiangular lines in $\mathbb{R}^n$ for $n=3 \cdot 2^{2t-1}, t\in{\mathbb N}.$
Thus, currently the best asymptotic results are summarized as follows:
  \begin{equation}\label{eq:ap}
     \frac29\le \limsup_{n\to\infty} \frac{M(n)}{n^2}\le \frac12,
  \end{equation}
where the upper bound is from \eqref{eq:gerzon}. The question of
the correct order of growth represents a difficult unresolved
problem. Contributing to the study of the asymptotic bounds, we
show that for $n=3(2k-1)^2-4$ and $\alpha=\frac{1}{2k-1},$ for all
integer $k\ge 2,$
  \begin{equation}\label{eq:b}
   M_\alpha(n)\le \frac{(n+1)(n+2)}{6}.
   \end{equation}
{\sc (Added in proof)} After this paper was accepted, C. Greaves et al. posted a pre\-print \cite{Greaves14} in which the upper bounds for $n=14,16$ were improved to $M(14)\leq 29$ and $M(16)\leq 41,$ respectively.

\remove{
\vspace*{.1in} In Section 2, we introduce
the definitions and properties of semidefinite programming.
Section 3 investigates the mathematics knowledge to derive matrix
inequalities in SDP. Section 4 will talk about the matlab
computation results running by CVX toolbox and we will list the
complete table to dimensions $n=97$. In section 5, we will talk
about tight spherical harmonic index $4$ design which are
equiangular line sets in $\mathbb{R}^n$, where $n$ is three times
odd number square minus four. We use linear programming method to
prove this type of equiangular lines sets in $\mathbb{R}^n$ such
that upper bounds for their size are $\frac{(n+1)(n+2)}{6}$.
However, we use SDP to prove that several lower dimension cases
are strictly less than this bound.}

\section{SDP bounds for equiangular lines}
Many problems in operations research, combinatorial optimization,
control theory, and discrete geometry can be modelled or approximated
as semidefinite programming. SDP optimization problems are usually stated in the
following form:
  \begin{equation*}
   \min  c^T x
  \end{equation*}
  \begin{equation*}
\text{\rm subject to } \quad F_0 + \sum_{i=1}^m F_i x_i  \succeq 0, \quad x\in\reals^m,
  \end{equation*}
where $c \in \mathbb{R}^m$ is a given vector of coefficients, $F_i,i=0,1,\dots$ are $n\times n$
symmetric matrices, and "$\succeq$" means that the matrix is positive semidefinite.
SDP problems fall in the class of convex optimization problems since the domain
of feasible solutions is a convex subset of $\reals^m.$ For the case
of diagonal matrices $F_i,$ SDP turns into a linear programming (LP) problem.
Properties of SDP problems and algorithms for their solution are discussed,
for instance, in \cite{boyd96}.
Most SDP solvers such as CSDP, Sedumi, SDPT3 use interior point methods
originating with Karmarkar's celebrated algorithm
(we used CVX toolbox in Matlab).

Let $\mathcal{C}\subset S^{n-1}$ be a set of unit vectors in $\reals^n$ such that $\langle x,x'\rangle\le a$
for all $x,x'\in \mathcal{C}, x\ne x'$ (a {\em spherical code}). As shown
by Bachoc and Vallentin \cite{bac08a}, the problem of estimating the maximum size of $\mathcal{C}$ can be stated as
an SDP problem. In particular, for $a=1/2,$ this is the famous ``kissing number problem'', i.e.,
the question about the maximum number of nonoverlapping unit spheres that can touch a given unit sphere.
A particular case of the main result in \cite{bac08a} was used in \cite{barg13} to find new bounds on
the maximum cardinality of spherical two-distance sets.

Let us introduce some notation.
Let $G_k^{(n)}(t), k=0,1,\dots$ denote the Gegenbauer polynomials
of degree $k$, i.e., a family of polynomials defined recursively as follows:
$G_0^{(n)}\equiv 1, G_1^{(n)}(t)=t,$
 and
$$
    G_k^{(n)}(t)=\frac{(2k+n-4)tG_{k-1}^{(n)}(t)-(k-1)G_{k-2}^{(n)}(t)}
   {k+n-3}, \quad k\ge 2.
$$
Following \cite{bac08a}, define a $(p-k+1)\times(p-k+1)$ matrix $Y_k^n(u,v,t), k\ge 0,$
  $$
    (Y_k^n(u,v,t))_{ij}=u^iv^j ((1-u^2)(1-v^2))^{k/2}
      G_k^{(n-1)}\Big(\frac{t-uv}{\sqrt{(1-u^2)(1-v^2)}}\Big),
  $$
where $p \in \mathbb{N}$, and a matrix $S_k^n(u,v,t)$ by setting
  \begin{equation*}
    S_k^n(u,v,t)=\frac16\sum_\sigma Y_k^n(\sigma(u,v,t)),
  \end{equation*}
where the sum is over all permutations of 3 elements. Note that
$(S_k^n(1,1,1))_{ij}=0$ for all $i,j$ and all $k\ge 1.$
Let $\mathcal{C}$ be a spherical code. As shown in \cite{del77b},
  \begin{equation}\label{eq:de}
\sum_{(x,y) \in \mathcal{C}^2}G_k^{(n)}(\langle x, y\rangle) \geq 0,
\end{equation}
and as shown in \cite{bac08a},
\begin{equation}\label{eq:be}
\sum_{(x,y,z) \in \mathcal{C}^3}S_k^n(\langle x, y\rangle, \langle x , z\rangle, \langle y, z\rangle) \succeq 0.
\end{equation}
Inequalities \eqref{eq:de} and \eqref{eq:be} can be used to formulate a general SDP problem for upper
bounds on the cardinality of spherical codes in $\reals^n$ \cite{bac08a}.

Using the approach of \cite{barg13}, we obtain the following SDP bound on $M(n).$
\begin{theorem}
Let $\mathcal{C}$ be set of a equiangular lines with inner product values
either $a$ or $-a$. Let $p$ be the positive integer. The
cardinality $|\mathcal{C}|$ is bounded above by the solution of
the following semi-definite programming problem :
\begin{equation}\label{eq:SDPF}
1+ \frac{1}{3} \max (x_1+x_2)
\end{equation}
subject to
     \begin{align}
&\Big( \begin{array}{@{\hspace*{.05in}}c@{\hspace*{.05in}}c@{\hspace*{.05in}}}
1 & 0  \\
0 &1
\end{array} \Big) +\frac{1}{3}
\Big( \begin{array}{@{\hspace*{.05in}}c@{\hspace*{.05in}}c@{\hspace*{.05in}}}
0 & 1  \\
1 & 1
\end{array} \Big) (x_1 +x_2) +
\Big( \begin{array}{@{\hspace*{.05in}}c@{\hspace*{.05in}}c@{\hspace*{.05in}}}
0 & 0  \\
0 & 1
\end{array} \Big) (x_3+x_4+x_5+x_6) \succeq 0  \label{eq:c1} \\
     & S^n_k(1,1,1) + S^n_k(a,a,1)x_1 + S^n_k(-a,-a,1)x_2 + S^n_k(a,a,a)x_3 \nonumber \\
              &  \hspace*{.8in}+S^n_k(a,a,-a)x_4+S^n_k(a,-a,-a)x_5
                    +S^n_k(-a,-a,-a)x_6 \succeq 0 \\
&\hspace*{2.5in}3+ G_k^{(n)}(a)x_1+ G_k^{(n)}(-a)x_2 \geq 0,\label{eq:c3}
     \end{align}
where $k=0,1, \cdots, p$ and $x_j \geq 0$, $j=1,\cdots,6$.
\end{theorem}
To compute bounds on $M(n)$, we found solutions of the SDP problem  \eqref{eq:SDPF}-\eqref{eq:c3},
restricting our calculation to the case $p=5.$  In Table \ref{table3} we list the values of SDP bounds for all possible
angles except the angle $\arccos\frac{1}{3}$ which is not included because of \eqref{eq:1/3} (note that
the SDP bounds for other angles are much greater than
$2(n-1)$). The column labelled `max' refers to the maximum of the SDP
bounds among all possible angles. The last column in the table gives the value of
the angle for which the maximum is attained.

{\footnotesize
\begin{table}
\begin{center}\begin{tabular}{|c|c|c|c|c|c|c|c|c|c|}
\hline $n$ &  1/5 & 1/7 & 1/9 & 1/11 & 1/13 &  1/15 & max & Gerzon &  angle\\
\hline

           22      &    176&       39       &29 &      26&      25&       24&          176&      253&       1/5\\
           23       &   276 &      42       &31  &    28  &     26&       25&          276 &         276 &         1/5\\
           24        &  276  &      46       &33  &    29  &    27&     26&          276    &      300    &       1/5\\
           25         & 276   &        50      & 35&        31&       29&           28       &   276       &   325       &     1/5\\
           26          &276    &   54       &37     &  32      & 30&      29&          276    &      351    &       1/5\\
           27 &         276     &  58           &40  &     34   &    31&      30&          276 &         378 &          1/5\\
           28  &        276      &     64      & 42   &    36    &   33&       31&          276 &         406 &           1/5\\
           29   &       276       &  69       &44      & 37       &  34&       33&          276  &        435  &          1/5\\
           30    &      276       &75      & 47        &39       &36&       34&          276      &    465      &     1/5\\
           31     &     276       &82        & 49       &41&       37&       35&          276      &    496      &     1/5\\
           32      &    276       &90       &52       &43&       39&       37&          276         & 528         &   1/5\\
           33       &   276        &   99           &55  &         45 &      40&         38      &   276          &561        &   1/5\\
           34        &  276        &108     &  57       &46&       42&      39&          276&          595          &  1/5\\
           35         & 276         & 120        &60     &  48&       43&       41&          276&          630       &     1/5\\
           36          &276       &132         &  64      & 50&       45&       42&          276 &         666        &   1/5\\
           37&          276        &  148     &  67       &52&       47&       44&          276   &       703          &  1/5\\
           38 &         276       &165       &70        &54&     48&      45&        276          &741            &1/5\\
           39  &        276        &187       &74       &57     &    50&       46&          276     &     780       &     1/5\\
           40   &       276       &213       &78       &59&       52&       48&          276        &  820          &  1/5\\
           41    &      276        &  246        &   82 &        61&       53&       49          &276    &      861 &         1/5\\
           42     &     276         & 288      & 86      & 63&       55&      51&       288   &       903           & 1/7\\
           43      &    276          &344       &90       &66&       57&       52&          344&          946        &    1/7\\
           44       &   276       & 422       &95       &68&       59&       54&        422&        990          &  1/7\\
           45        &  276        &  540     &     100 &      71&       60&           56   &       540         &1035     &       1/7\\
           46         & 276         & 736   &    105     &    73&       62&       57&          736&         1081 &          1/7\\
           47          &276    &     1128  &     110      & 76&       64&       59&         1128   &      1128    &        1/7\\
           48&          276     &    1128       &116       &78&       66&       60&         1128    &     1176     &       1/7\\
           49 &         276      &   1128        &122       &81&         68&       62&         1128  &       1225   &         1/7\\
           50  &        276       &  1128       &129       &84&       70&           64&         1128  &       1275   &         1/7\\
           51   &       276        & 1128        &  136     &  87&        72&       65&         1128   &      1326    &        1/7\\
           52    &      276         &1128       &143       &90&       74&      67&       1128         &1378            &1/7\\
           53     &     276&         1128       &151       &93&      76&       69&         1128        & 1431           & 1/7\\
           54      &    276 &        1128       &   160     &  96&       78&      70&         1128      &   1485         &   1/7\\
           55       &   276  &       1128       &169         & 100&       81&       72&         1128     &    1540        &    1/7\\
           56        &  276   &      1128        &179       &103&       83&       74&         1128        & 1596           & 1/7\\
           57         & 276    &     1128        &  190      & 106&         85&           76&         1128 &       1653    &        1/7\\
           58          &276     &    1128       &201       &110   &   87&     77&         1128        & 1711     &      1/7\\
           59&          276      &   1128       &214       &114&       90&       79&         1128       &  1770      &      1/7\\
           60 &         276       &  1128      & 228       &118&    92&       81&         1128         &1830          &  1/7\\
           61  &       279         &1128         & 244      &    122&       94&       83         &1128         &1891    &       1/7\\
           62   &       290         &1128      & 261       & 126&       97&      85        &1128       &  1953           & 1/7\\
           63    &   301&         1128        &  280       &130&       99&       87         &1128       &  2016           & 1/7\\
           64     &  313 &        1128       &301       &134&    102 &      89         &1128        & 2080            &1/7\\
           65      &    326 &        1128        &  325  &     139&          105           &91         &1128        & 2145&            1/7\\
           66       &339   &      1128          &352      &    144&       107       &92         &1128    &     2211 &          1/7\\
           67        &  353  &       1128       &382       &148&       110       &94         &1128 &        2278        &    1/7\\
           68       &367      &   1128       &418      &153&       113&       97&         1128     &    2346           & 1/7\\
           69       &382       &  1128       &   460      & 159&       115&       99&         1128&         2415          &  1/7\\
           70       &398        & 1128      & 509       &164&       118&       101&        1128   &      2485           & 1/7\\
           71&          416      &   1128      &    568  &      170&       121&     103&        1128 &        2556       &     1/7\\
           72 &         434       &  1128     &     640   &    176&        124&       105&         1128 &        2628       &     1/7\\
           73  &     453&         1128       &   730       & 182&       127&      107&         1128  &       2701           & 1/7\\
           74   &    473 &        1128      & 845       &188&      130&      109&         1128        & 2775          &  1/7\\
           75    &   494  &       1128     &    1000     &  195&       134&          112&         1128  &       2850   &         1/7\\
76     &  517   &      1128    &     1216      & 202&       137&       114&         1216 &        2926        &    1/9\\
           77      &    542 &        1128      &   1540    &      210&       140    &   116&         1540  &         3003  & 1/9\\
           78       &568     &    1128        & 2080       &217&         144    &   118&        2080      &   3081       &     1/9\\
           79        &  596   &      1128    &     3160     &  225&       147&      121&        3160      &   3160          & 1/9\\
           80         & 626    &     1128   &      3160      & 234&       151&       123&         3160      &   3240    &        1/9\\
           81       &658        & 1128     &    3160          &243&       154&          126      &   3160  &       3321 &           1/9\\
           82        &693        & 1128   &      3160      & 252&      158&      128&         3160        & 3403         &  1/9\\

\hline
\end{tabular}\end{center}
\end{table}
\begin{table}
\begin{center}\begin{tabular}{|c|c|c|c|c|c|c|c|c|c|}
\hline $n$ &  1/5 & 1/7 & 1/9 & 1/11 & 1/13 &  1/15 & max & Gerzon &  angle\\
\hline
                      82        &693        & 1128   &      3160      & 252&      158&      128&         3160        & 3403         &  1/9\\
           83          &731       &  1128&         3160     &  262&       162&       130&         3160&         3486      &      1/9\\
           84&       772    &     1128  &       3160       &272&       166&      133&         3160 &        3570           & 1/9\\
           85 &      816     &    1128         &3160      & 283&          170&          136&         3160   &      3655       &     1/9\\
           86  &        866   &     1128      &   3160     &  294&       174&       138&         3160      &   3741            &1/9\\
           87   &       920    &     1128    &     3160     &  307&       178&       141&         3160&         3828&           1/9\\
           88    &   979        & 1128      &   3160        &  320 &      182&       143&        3160  &       3916  &         1/9\\
           89     &    1046      &   1128  &       3160     &  333&        186&       146&         3160 &        4005 &          1/9\\
           90      & 1120         &1128   &      3160       &348&     191&    149&         3160&         4095         &   1/9\\
           91       &1203&         1128  &       3160       &   364 &         196&       152&         3160&         4186&           1/9\\
           92        & 1298&         1128        & 3160     &  380&       200&     154&        3160        & 4278        &    1/9\\
           93         &1406 &        1128       &  3160     &  398&       205&       157&         3160      &   4371      &      1/9\\
           94&       1515    &     1128        & 3160       &417&       210&       160&         3160         &4465         &   1/9\\
           95 &      1556     &    1128       &  3160       &438&       215&       163&         3160         &4560          &  1/9\\
           96  &     1599      &   1128      &   3160       & 460&      220&        166&        3160      &   4656            &1/9\\
           97   &    1644       &  1128     &    3160       &   485&       226&       169&         3160      &   4753         &   1/9\\
           98    &   1691        & 1128    &     3160       & 511&      231&      172&      3160     &    4851           & 1/9\\
           99     &  1739         &1128   &      3160       &   540&        237&          176&         3160  &         4950 &          1/9\\
          100      & 1790  &       1128  &       3160       &571&      243&       179&       3160      &   5050          &  1/9\\
          101       &1842   &      1128 &        3160       &   606&       249&       182&         3160  &       5151&            1/9\\
          102&       1897    &     1128&         3160       &644&       255&       185&         3160     &    5253     &       1/9\\
          103 &      1954     &    1128         &3160       &686&       262&       189&         3160      &   5356 &           1/9\\
          104  &     2014      &   1128        & 3160       &734&       268&      192&       3160       &  5460      &      1/9\\
          105   &      2077     &    1128     &    3160     &   787&       275&          196&         3160 &        5565 &           1/9\\
          106    &   2142        &1128       &  3160        &  848 &      282&       199&         3160       &  5671       &    1/9\\
          107     &    2211&         1128   &      3160     &  917 &     289&       203&        3160         &5778&            1/9\\
          108      & 2282   &      1128    &     3160       &997   &    297&       206&         3160     &    5886 &           1/9\\
          109       &2358    &     1128   &      3160       &  1090  &      305&       210&         3160     &    5995&           1/9\\
          110 &      2437     &    1128  &       3160       &  1200 &      313&       214&         3160      &   6105   &         1/9\\
          111  &     2521      &   1128 &        3160        & 1332&       321&       218&         3160       &  6216 &            1/9\\
          112   &    2609       &  1128&         3160       &1493&      330&       222&     3160         &6328          &  1/9\\
          113    &   2702        & 1128         &3160       &  1695&          339&          226&         3160 &        6441  &         1/9\\
          114     &  2800         &1128        & 3160       &1954     & 348&      230&       3160       &  6555            &1/9\\
          115      & 2904&         1128       &  3160       &  2300     &  357&       234&        3160   &      6670 &           1/9\\
          116       &  3015&         1128    &     3160     &    2784  &      367&       238&         3160       &  6786 &           1/9\\
          117      & 3132   &      1128     &    3160       &  3510   &       378 &      242&        3510         &6903  &1/11 \\
          118       &  3257  &       1128  &       3160     &    4720&       388&       247&        4720        & 7021   &        1/11 \\
          119       &  3390   &      1128 &        3160     &  7140&   399&     251&       7140&         7140    &       1/11 \\
          120      & 3532      &   1128  &       3160       &  7140      & 411&          256&         7140  &       7260&          1/11 \\
          121       &3684       &  1128 &        3160       &  7140     &   423&       260&         7140      &   7381    &       1/11 \\
          122        & 3848      &   1128       &  3160     &    7140  &     436&       265&         7140      &   7503&            1/11 \\
          123        & 4024       &  1128      &   3160     &    7140 &      449&       270&         7140        & 7626  &         1/11 \\
          124        & 4214        & 1128     &    3160     &    7140&       462&       275&         7140         &7750&           1/11 \\
          125       &4419         &1128      &   3160       &  7140&       477&          280&         7140      &   7875 &           1/11 \\
          126  &     4643&         1128     &    3160       &  7140&       492&       285&         7140        & 8001  &          1/11 \\
          127   &      4887 &        1128  &       3160     &  7140&         508 &      290&       7140 &         8128  &        1/11 \\
          128    &   5153  &       1128   &      3160       &  7140      & 524&      295&    7140        & 8256          & 1/11 \\
          129     &    5447  &       1128&         3160     &    7140   &     541&          301 &        7140 &        8385&           1/11 \\
          130     &  5770     &    1128 &        3160       &  7140     &     560 &      306&         7140       &  8515     &       1/11 \\
          131   &    6130      &   1128&         3160       &  7140     &  579&       312&         7140        & 8646  &          1/11 \\
          132    &   6531       &1130 &        3160         &7140       &599&      317&        7140   &      8778        &   1/11 \\
          133     &  6982       &1158         &3160         &7140       &620&       323&         7140  &       8911  &         1/11 \\
          134      & 7493&       1187        & 3160         &7140      &  643&       329&     7493&         9045     &       1/5\\
          135       &8075 &        1218     &    3160       &  7140   &    667&          336&       8075&         9180&            1/5 \\
          136&       8747  &     1249      &   3160         &7140    &   692&      342&       8747&         9316       &     1/5\\
         *137 &      9528   &      1282   &      3160       &  7140 &      719&       348&       9528&         9453     &       1/5 \\
          *138  &      10450  &     1315  &       3160       &  7140&       747&       355&        10450&         9591    &        1/5 \\
          *139   &     11553   &    1350 &        3160       &7140 &       778&       362&        11553  &       9730      &      1/5 \\

 \hline
\end{tabular}
\end{center}\caption{Values of the SDP bound on $M(n), 22\le n\le 139$}
\label{table3}
\end{table}}

Some comments on the tables are in order.
Observe that $M_{1/5}=276$ for $ 23 \leq n \leq
60$. Combined with the results of \cite{lem73}, this implies that
$M(n)=276$ for $23 \leq n \leq 41$ and $M(43)=344$. The case
$n=42$ remains open since we only obtain that $276 \leq M(42) \leq 288$ for the angle
$\arccos 1/7$.

Improvements of the Gerzon upper bound \eqref{eq:gerzon} are obtained for $n \leq 136.$
The last 3 entries in Table \ref{table3} produced no improvements, and are marked by an asterisk because of that.
Similarly, the SDP problem yielded no improvements for higher dimensions.

An interesting, unexplained observation regarding this table is
that the SDP bound for $M_\alpha(n)$ has long stable ranges for
dimensions starting with the value $n= d^2 -2,$ where $d$ is an
odd integer and $\alpha=1/d.$ For instance, one such region begins
with $d=5,$ another with $d=7.$ The same phenomenon can
observed for $d=9$  where the SDP value $M_\alpha(n)\le 3160$ is
obtained for all values of $n$ satisfying $ 79\le n\le 227$ and
for $d=11$ where the value $7140$ appears for all $n, 119\le n\le
347$.

Note that the SDP bound gives the same value as the Gerzon bound
for $n=47, 79$ and $119,$ and that these three dimensions are of
the form $n=(2k-1)^2-2,$ where $k \geq 2$ is a positive integer.
Bannai, Munemasa, and Venkov
\cite{ban04} showed that for $n=47,79$ the maximum possible size
$M(n)$ cannot attain this value while the case $n=119$ is still
open. The result of \cite{ban04} relies on the fact that an equiangular
line set in $\mathbb{R}^n$ with cardinality $\frac{n(n+1)}{2}$
gives rise to a spherical two-distance set of size $(n-1)(n+2)/2$
in $\mathbb{R}^{n-1},$ and such sets are related to tight
spherical $4$-designs whose existence can be sometimes ruled out.

Based on the earlier results and our calculations, we make the following

\vspace*{.05in}{\sc Conjecture:} 
{\em There exist $1128$ equiangular lines in $\mathbb{R}^{48}$ with
angle $ arccos(1/7)$ and $3160$ equiangular lines in
$\mathbb{R}^{80}$ with angle $ arccos(1/9)$.}

\vspace*{.05in}If this conjecture is true, then $M(n)=1128$ for $ 48 \leq n \leq 75$
and $M(n)=3160$ for $ 80 \leq n \leq 116$.

\remove{$M(n)$ has been twice to have stable values for several dimensions
starting from $n=7$ and $23$. Namely, $M(n)=28$ for $ 7 \leq n
\leq 13$ and $M(n) =276$ for $23 \leq n \leq 41$. We expect
similar phenomenon for higher dimensions. If Conjecture 2.2 is
true, then $M(n)=1128$ for $ 48 \leq n \leq 75$ and $M(n)=3160$
for $ 80 \leq n \leq 116$. We also believe that this is general
phenomenon for all $n$ when $n$ is around odd integer square.}

\section{Tight spherical designs of harmonic index $4$ and equiangular lines}

\begin{definition} Let $t$ be a natural number.
A finite subset $X$ of the unit sphere $S^{n-1}$
is called a \emph{spherical $t$-design} if, for any polynomial
$f(x)=f(x_1,x_2, \ldots, x_n)$ of degree at most $t$, the
following equality holds :
$$
\frac{1}{|S^{n-1}|} \int_{S^{n-1}}f(x) d\sigma(x) =
\frac{1}{|X|}\sum_{x \in X} f(x).
$$
\end{definition}
A spherical $t$-design is called tight if it attains the LP bound of \cite{del77b},
also called the {\em absolute bound}.

An equivalent definition of spherical designs can be given in terms of harmonic polynomials.
Let
$\text{Harm}_t(\mathbb{R}^n)$ be the set of homogeneous harmonic
polynomials of degree $t$ on $\mathbb{R}^n$. Then the set $X$ is a spherical design \cite{del77b} if
 $$
  \sum_{x \in X} f(x) = 0 \quad\forall  f(x) \in \text{Harm}_j(\mathbb{R}^n), 1 \leq j \leq t .
$$

The following definition was recently proposed by Bannai, Okuda,
and Tagami  \cite{ban13}: A {\em spherical design
of harmonic index $t$} is a finite subset $X\subset S^{n-1}$ such that
  $$
    \sum_{x \in X} f(x) = 0 \quad\forall  f(x) \in \text{Harm}_t(\mathbb{R}^n).
$$

An LP bound for spherical designs of harmonic index $t$ was derived in
 \cite{ban13}. Similarly, if this bound is attained, then the design is called tight.
 Our interest in tight spherical designs of a fixed harmonic index is motivated by a result in \cite{ban13}
 which shows that a tight design of index $4$ gives rise to an equiangular line set in $\reals^n$
 with angle $\arccos a=\sqrt{3/(n+4)}$. Since $a=\frac{1}{2k-1}$ for some integer
$k \geq 2,$ we find that $n=3(2k-1)^2-4$. These considerations motivate the following result.
\begin{theorem} Let $n=3(2k-1)^2-4, k\ge 2.$ The cardinality $N$ of any equiangular line set in $\reals^n$
  with inner product $a=1/(2k-1)$ satisfies the inequality
\begin{equation}\label{eq:S1}
|S| \leq \frac{(n+1)(n+2)}{6}.
\end{equation}
\end{theorem}
\proof

To prove this result we use the LP bound of \cite{del77b} that has the following form:
{\em Let $T\subset[-1,1]$. Let $S=\{x_1,x_2,\dots,x_N\}$ be a set of unit
vectors in $\mathbb{R}^n$ such $\langle x_i,x_j\rangle\in T\cup\{1\}.$
Let $f(t)=\sum_k f_kG^n_k(t)$ be a polynomial such that $f_0>0, f_k\ge 0, k\ge 1$
and that $f(t) \leq 0$ for all $t \in T.$ Then}
  \begin{equation}\label{eq:LP}
   |S| \leq \Big\lfloor \frac{f(1)}{f_0} \Big\rfloor.
  \end{equation}

Consider the polynomial
   $$
   f(t) = (t^2-a^2)\Big(t^2 +\frac{a^2n+4a^2-6}{n+4}\Big).
   $$
Let $X\subset \reals^n$ be an equiangular line set with inner product $a$. Then $T=\{\pm a\},$ and $f(t)=0$ for $t\in T.$
Computing the Gegenbauer expansion of $f(t),$ we obtain
     \begin{align*}
&f_0 = -\frac{a^4n^2 + 6a^2n(a^2-1) + 8a^4 - 6a^2(n+2) +
3}{n^2 + 6n + 8}\\
&f_1 = f_2 = f_3 =0\\
&f_4 =\frac{n^2 - 1}{(n+ 2)(n + 4)}
    \end{align*}
We need to check that $f_0>0.$ Substituting the values of $n$ and $a$, we obtain
   $$
   f_0 =\frac{8k(k -1)}{(2k - 1)^4(12k^2 - 12k + 1)}\ge 0 \quad \text{for }k\ge 2.
   $$
Thus, $f(t)$ satisfies the conditions of the LP bound, and we obtain
   $$
   |S| \leq \frac{f(1)}{f_0} = \frac{(a^2 - 1)(n + 2)(n + a^2n + 4a^2
- 2)}{a^4n^2 + 6a^4n + 8a^4 - 6a^2n - 12a^2 + 3}.
   $$
In particular, putting $a=\frac{1}{2k-1}$ and $n=3(2k-1)^2-4 =
12k^2-12k-1$, we obtain
  $$
  \frac{f(1)}{f_0}=\frac{(n+1)(n+2)}{6}. \text{\qed}
   $$

This theorem gives
infinitely many values of $n$ for which the upper bound $M_a(n)$ is strictly
less than the Gerzon bound, yielding the asymptotic constant 1/6 for the growth rate of the
quantity $M_\alpha(n)$ (cf. \eqref{eq:ap}-\eqref{eq:b}).

\begin{remark} Observe that the relative bound \eqref{eq:relative} is an instance of the LP bound
\eqref{eq:LP};
see \cite{del77b}. Thus, the SDP bound \eqref{eq:SDPF}-\eqref{eq:c3} is as strong or stronger than the
bound  \eqref{eq:relative}.
\end{remark}

\begin{remark}
Using SDP, we further show that for some dimensions the LP bound \eqref{eq:S1} cannot be attained.
Indeed, for $k=3,4,5$ we obtain the values of the dimension $n=71,143,239,$ respectively,
and the SDP bound implies that
  $$
     M_{1/5}(71)\le 416, \quad M_{1/7}(143)\le 1506, \quad M_{1/9}(239)\le 3902,
  $$
  which is much smaller than the values $876,3480,9640$ obtained from \eqref{eq:S1}.
Extending these calculations, we have shown that for $k\le 54$ and $n=3(2k-1)^2-4\le 34343$ the
SDP bound improves upon the LP bound \eqref{eq:S1}.
\end{remark}

In conclusion, we note that the value of the maximum in the LP problem for the maximum cardinality
of equiangular line sets with a given angle can be explicitly characterized.  The LP problem has the
following form:
  \begin{equation}\label{eq:max}
   M_a(n)\le \max\{1+x_1+x_2, x_1\ge 0, x_2\ge 0\}
  \end{equation}
   subject to
\begin{align}\label{eq:co}
 1+G^n_k(a)x_1+G^n_k(-a)x_2 \geq 0 \quad \text{for $k= 1, 2, \dots .$}
\end{align}
\begin{theorem}
 Let $a\in(0,1),$
     \begin{equation}\label{eq:gn}
     g_n = \min_{k\ge 0} \frac{1}{|G^n_k(a)|}
     \end{equation}
where $k$ is even and such that $G^n_k(a)<0$.
Then
     $$M_a(n) \leq g_n +1,
     $$
     where the value $g_n+1$ is the solution of the LP problem \eqref{eq:max},\eqref{eq:co}.
\end{theorem}

\proof
Let $k$ be even, then $G_k^n(t)$ is an even function, so inequalities \eqref{eq:co} take the form
   \begin{equation}\label{eq:in}
   1+ G^n_k(a)(x_1+x_2) \geq 0, \quad k=2m, m\in{\mathbb N}.
   \end{equation}
These inequalities define a set of half-planes whose boundaries are parallel to the objective
function. The inequalities for odd $k$ are bounded by lines that are perpendicular to the boundaries of
the even-indexed constraints, and therefore can be disregarded.
We conclude that the maximum is attained on the line $1+G^n_k(a)(x_1 + x_2)=0$ for some even $k.$
The inequalities with $k$ such that $G^n_k(a) \geq 0$ are trivially satisfied, therefore,
we consider only those values of $k$ when
$G^n_k(a) < 0$. Eq.~\eqref{eq:in} implies that, for all even $k$,
\begin{align*}
x_1+x_2 \leq -\frac{1}{G^n_k(a)} = \frac{1}{|G^n_k(a)|}.
\end{align*}
This completes the proof.
\qed

\vspace*{.1in} To give an example of using this theorem, take
$n=71$ and $a=\frac{1}{5}.$ To find a bound on $M_a(n),$ we
estimate the quantity $g_n$ in \eqref{eq:gn} by computing
   $$
   \min_{0\le k\le 100} \frac 1{|G_k^{(71)}(1/5)|}
   $$
for all even $k$ such that $G_k^{(71)}(1/5)<0.$ The smallest value
is obtained for $k=4,$ and $G^{(71)}_4(1/5)=-1/875.$ Thus, we
obtain $M_{1/5}(71)\le 876$. Of course, it could be possible that
for greater $k$ we obtain a smaller value of the bound, but this
is not supported by our experiments (although we do not have a
proof that $k=4$ is the optimal choice).

Experiments also suggest that $k=4$ may be the universal
optimal choice for infinitely many values of $n$ and
$a$. Indeed, we have
  $$
    G^n_4(x)=\frac{(n+2)(n+4)x^4-6(n+2)x^2+3}{n^2-1}.
    $$
Taking $n=3(2t-1)^2-4$ and $a=1/(2t-1),$ where $t
\geq 2,$  we obtain the expression
   $$
\frac{1}{G^n_4(a)}+1 = 2t(t-1)(12t^2-12t+1)=
\frac{(n+1)(n+2)}{6}
  $$
which coincides with the LP bound \eqref{eq:S1}.

\vspace*{.1in}{\sc Acknowledgement:} We thank Peter Casazza, Eiichi Bannai, and Ferenc
Sz\"{o}ll\H{o}si for insightful discussions on the topics of this paper.


Alexander Barg was supported in part by NSF grants DMS1101697, CCF1217245, and CCF1217894, and by NSA grant H98230-12-1-0260. Wei-Hsuan Yu was supported in part by NSF grant CCF CCF1217245.

\providecommand{\bysame}{\leavevmode\hbox
to3em{\hrulefill}\thinspace} \providecommand{\href}[2]{#2}
\bibliographystyle{amsalpha}

\end{document}